\pdfoutput=1 
\documentclass[letterpaper,conference,10pt]{IWMRRF}
\usepackage[pdftex]{graphicx}
\usepackage{graphicx}  

\usepackage{url}       

\usepackage{amsmath,amssymb}   
\usepackage{array}

\usepackage{pgfplots}
 \usetikzlibrary{arrows.meta}
\pgfplotsset{compat=newest}
\usetikzlibrary{external}
\usepackage{color}
\usepackage{commath}
\usepackage{cite}

\newtheorem{theorem}{Theorem}

\newtheorem{definition}{Definition}

\begin{document}

\title{Characterization of Separatrices in Holomorphic Dynamical Systems}

\author{\authorblockN{Marcus Heitel\authorrefmark{1}, Dirk Lebiedz\authorrefmark{1}}
\authorblockA{\authorrefmark{1}Institute for Numerical Mathematics, Ulm University, Germany}}

\maketitle

%
\begin{abstract}
	Multiple time scales in dynamical systems lead to a bundling of trajectories onto slow invariant manifolds (SIMs). Although they are absent in two-dimensional holomorphic dynamical systems, a bundling of orbits is often observed as well. They bundle onto special trajectories called separatrices. We apply numerical methods for the approximation of SIMs to holomorphic flows and show how a separatrix between two regions of periodic orbits can be characterized topologically. Complex time reveals a new perspective on holomorphic dynamical systems.
\end{abstract}

\section{Introduction}
Separatrices are structures similar to slow manifolds that characterize the flow of a holomorphic function. Their identification has become important in the context of the Riemann Hypothesis (RH) (see \cite{BroughanStructure,BroughanZeta,Schleich}). The challenging task for a possible proof of the RH is to characterize separatrices for the 1D-holomorphic flow of the Riemann zeta function \cite{SchleichEq}. 

One-dimensional holomorphic dynamical systems can be interpreted as dynamical systems in $\mathbb{R}^2$ with a special structure given by the Cauchy-Riemann equations. 
If a holomorphic function $f:\mathbb{C} \rightarrow \mathbb{C},$ $f:(x+iy) \mapsto u(x,y) + iv(x,y)$ is regarded as two-dimensional function $\hat{f}:\mathbb{R}^2 \rightarrow \mathbb{R}^2$ that maps $(x,y)^T$ to 
$\left(u(x,y),v(x,y)\right)^T$, then the Cauchy-Riemann equations yield
\begin{align}
	\hat{f}'(x,y) = \left(\begin{array}{rr} u_x & -v_x\\v_x &u_x 	\end{array}\right).
\end{align}

\section{Separatrix and the $\cosh$ flow}
In this holomorphic context separatrices can be defined as follows.
\begin{definition}[Broughan, \cite{BroughanStructure}] A trajectory $\gamma$ is a {\em positive (negative) separatrix} if for some $z\in \gamma$ the maximum interval of existence of the path commencing at $z$ and proceeding in positive (negative) time is finite. A trajectory $\gamma$ is a {\em separatrix} if it is a positive or negative separatrix.
\end{definition}

With this definition Broughan \cite{BroughanStructure} argues that separatrices are boundary components of special regions. These regions are the union of all trajectories with the qualitatively same phase space behavior, e.g., periodic orbits around the same center.

\begin{figure}[tp]
	\centering
	\includegraphics[scale=0.95]{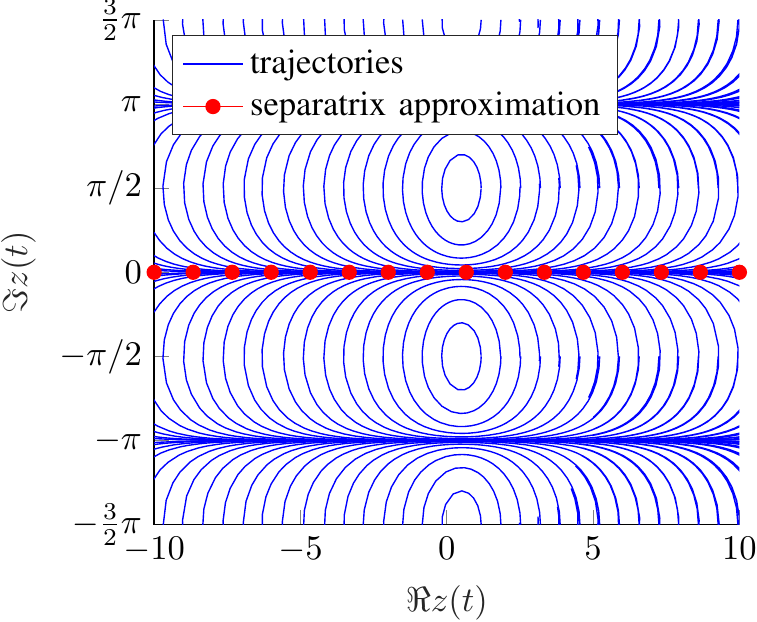}
	\caption{Phase portrait for the holomorphic flow $\dot{z} = \cosh\left(z-\frac{1}{2}\right)$ for $z=a+b\,i$ with $a \in [-10,10] $ and $b \in \left[-\frac{3}{2}\pi, \frac{3}{2}\pi\right]$.}
	\label{fig:cosh}
\end{figure}

In Figure \ref{fig:cosh} trajectories of
\begin{equation}\label{eq:cosh}
	\dot{z}=\cosh(z-0.5)
\end{equation} 
bundle near the separatrices which are the lines $\Im z = k\pi$ for every $k\in \mathbb{Z}$. This phase space behavior of trajectories is similar to those of slow-fast systems. The difference is here that separatrices are in general not slow (invariant) manifolds. Indeed, there is no spectral gap because of the Cauchy-Riemann equations.

\section{Application of SIM Methods}
Due to the common bundling behavior, we apply slow manifold approximation methods like the zero-derivative-principle (ZDP) \cite{Zagaris2009} and the following boundary value problem of Lebiedz and Unger \cite{Lebiedz2016} to holomorphic dynamical systems:
\begin{subequations} \label{formula:Unger}
	\vspace*{-1em}
	\begin{alignat}{2}
	\min\limits_{z(\cdot)}    & \quad & 				& \hspace*{-5mm} \left\|\ddot{z}(t_0)\right\|_2^2 \\
	\text{ s.t. } 							&  & \dot z 	  &= f(z), \quad t \in [t_0,t_1] \\
	&  &  x^* &= \Re z(t_1).
	\end{alignat}
\end{subequations}
The latter calculates pointwise a separatrix candidate for fixed real part of $z$ at time $t_1$, which corresponds to the reaction progress variable in the calculation of slow manifolds. The numerical results of this boundary value problem for several values of $x^* $ are plotted in Figure \ref{fig:cosh}. They all lie exactly on the separatrix. The same holds for the ZDP of order 1 applied to the imaginary part as variable, i.e., 
\begin{equation}
	\frac{\dif}{\dif t}\Im z = 0.
\end{equation}
In the special case 
\begin{align}
	f(z) &= \cosh(x-0.5+iy)  \\
	&= \cosh(x-0.5)\cos(y) + i \sinh(x-0.5)\sin(y), \nonumber
\end{align}
the ZDP yields the curves $\Im z = k\cdot \pi$ as well as $\Re z=0.5$.
However, in general these methods only lead to approximations of separatrices.  

\section{Separatrix Characterization}
\begin{figure}[tp]
	\centering
	\includegraphics[scale=0.95]{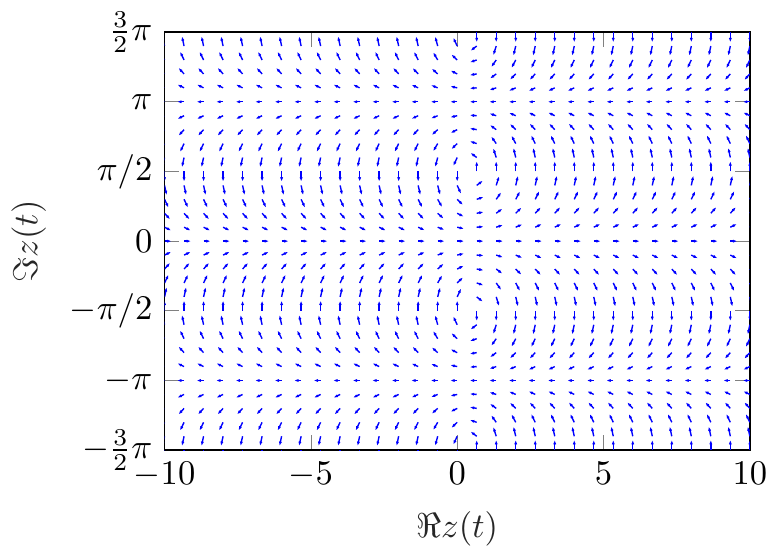}
	\caption{Direction field of $f(z) = \cosh(z-0.5)$}
	\label{fig:quiver}
\end{figure}

The idea is to characterize a separatrix by a closer look on the topology of the phase space near separatrices. Figure \ref{fig:quiver} depicts the direction field of  (\ref{eq:cosh}). 
Periodic orbits above and below the separatrices have different mathematical orientation. More precisely, the index of those curves differs, i.e., the number of (mathematically positive) rotations of a tangent vector along the curve. The following theorem captures this feature to characterize separatrices between two regions of periodic orbits around centers.

\begin{theorem}
	Let $z_0,z_1$ be two (neighbored) centers of the flow
	\begin{equation}\label{eq:ODE}
	\dot{z}(t)=f\!\left(z(t)\right)
	\end{equation} with a holomorphic function $f$. Let the regions 
	\begin{align}
	P_i := &\{z_i\} \cup\left\{z^*\in \mathbb{C} \,: \text{ solution of (\ref{eq:ODE}) with } z(0)=z^*\right. \nonumber\\
	& \left. \text{ is a periodic orbit around }z_i \,\right\}, \,\, i\in \left\{0,1\right\}
	\end{align}
	have common boundary component (i.e. separatrix, cf. Broughan \cite{BroughanStructure}) $\emptyset \neq C = \partial P_0 \cap \partial P_1 \cap U$ for an open set $U$. Then it holds  $z^* \in C \nonumber$
	\begin{align}
	\Leftrightarrow & \text{for }\varepsilon > 0 \text{ small enough } \exists \, z_i^* \in U_\varepsilon(z^*) \cap P_i,  \;\,  i\in \{0,1\} \nonumber\\
	& \text{ such that the orbits } \gamma_i \text{ of the flow (\ref{eq:ODE})} \text{ with } z(0)=z_i^*  \nonumber\\
	& \text{ satisfy } \,\, \text{ind}_{\gamma_0}(z_0) \cdot \text{ind}_{\gamma_1}(z_1) = -1.
	\end{align}
\end{theorem}

In case of system (\ref{eq:cosh}) it is clear that the indices of all periodic solution orbits around the same center are constant. Thus it changes only at the the curves $\Im z = k \cdot \pi$ for integers $k$, which are separatrices. 

\section{Complex Time}
In Equation (\ref{eq:ODE}) $f$ is a complex-valued function $f$ with complex argument. It seems natural to consider the solution $z(\cdot)$ of (\ref{eq:ODE}) as a function over a complex variable, i.e., complex time. If the real time is replaced by $t=\tau_1 + i \tau_2$ with $\tau_1,\tau_2 \in \mathbb{R}$, then the solution of 
\begin{equation}
	\frac{\dif}{\dif t}z(t) = f\!\left(z(t)\right)
\end{equation}
assuming $f$ to be complex differentiable w.r.t. time, satisfies
\begin{align}
	\frac{\partial}{\partial \tau_1}z(t) = f\!\left(z(t)\right) \;\text{ and }\; \frac{\partial}{\partial \tau_2}z(t) = i f\!\left(z(t)\right).
\end{align}
$\tau_1$ is called {\em real time} and $\tau_2$ {\em imaginary time}. Since multiplication by $i$ in $\mathbb{C}$ can be interpreted as rotation by $\frac{\pi}{2}$ in  $\mathbb{R}^2$, real time trajectories and imaginary time trajectories intersect orthogonally. Imaginary time allows to traverse separatrices of the real time flow and vice versa. 

The holomorphic derivative $f'(z_0)$ of a simple zero $z_0$ of $f$ determines the type of the equilibrium. If $f'(z_0)$ is real resp. imaginary, then it is a node resp. center. Otherwise, it is a focus. A node of the real time flow becomes a center of imaginary time flow and vice versa. Thus, switching the sign of $f$ and proceeding to imaginary time provide ways of ``manipulating'' the stability and attractivity of the dynamical system in the neighborhood of a simple zero. 

Separatrices of the real time flow can be approximated by maximizing the curvature of the pure imaginary time trajectories. For each (imaginary time) trajectory this yields a point close to the separatrix. 

\section*{Acknowledgment}
Special thanks goes to the Klaus-Tschira foundation for financial funding of the project.


\end{document}